\documentclass[times,12pt]{article}

\def\int{\mathbb{Z}}

\def\OO{{\mathcal O}}

\def\proof{{\bf Proof. }}

\def\Pf{\proof} 
 
\def\pf{\proof}

\usepackage{times}
\usepackage{graphics}
\usepackage{amssymb}
\usepackage{amscd}
\usepackage{amsmath}
\input xy 
\xyoption{all}
\title{On sheets of  conjugacy classes in good characteristic}
\newtheorem{theorem}{Theorem}[section]
\newtheorem{lemma}[theorem]{Lemma}
\newtheorem{corollary}[theorem]{Corollary}
\newtheorem{proposition}[theorem]{Proposition}
\newtheorem{definition}[theorem]{Definition}
\newtheorem{remark}[theorem]{Remark}
\newtheorem{example}[theorem]{Example}

\author{Giovanna Carnovale, Francesco Esposito\\
Dipartimento di Matematica Pura ed Applicata\\
Torre Archimede - via Trieste 63 - 35121 Padova - Italy\\
email: carnoval@math.unipd.it, esposito@math.unipd.it }

\date{}

\begin{document}
\maketitle
\begin{abstract}
We show that the sheets for  a connected reductive algebraic group $G$ over an algebraically closed field in good characteristic acting on itself by
conjugation are in bijection with $G$-conjugacy classes of triples $(M, Z(M)^\circ t, \OO)$ where $M$ is the connected centralizer of a semisimple element in $G$, 
$Z(M)^\circ t$ is a suitable coset in $Z(M)/Z(M)^\circ$ and $\OO$ is a rigid unipotent conjugacy class in $M$.
Any semisimple element is contained in a unique sheet $S$ and $S$ corresponds to a triple with $\OO=\{1\}$.
The sheets in $G$ containing a unipotent conjugacy class are precisely those corresponding to triples for which $M$ is a Levi subgroup of a parabolic subgroup of $G$ and such a  
class is unique.
\end{abstract}

\section{Introduction}

Given a regular action of  an algebraic group on a variety, the analysis of its orbits is a key point for the comprehension of the action. In order to understand the orbits one may want to put together those orbits sharing the same nature. One way to do that is to consider sheets. These are the irreducible components of the union of orbits of fixed dimension. 

The sheets for the adjoint action of a complex connected reductive group $G$ on its Lie algebra ${\mathfrak g}$ have been extensively studied. Those containing a semisimple element are called Dixmier sheets and they appeared at first in \cite{dix}.  A thorough analysis of sheets in a complex reductive Lie algebra is to be found in \cite{bo, BK,katsylo}. The sheets in this case are explicitly described, they are in bijection with $G$-orbits of  pairs $({\mathfrak l}, \OO)$ where ${\mathfrak l}$ is a Levi subalgebra of ${\mathfrak g}$ and $\OO$ is a nilpotent orbit in $[{\mathfrak l},\,{\mathfrak l}]$ which is itself a sheet (rigid nilpotent orbit).  This parametrization of sheets is  naturally related to the notion of induced orbits introduced in \cite{lusp}.  Semisimple elements lie in a unique sheet and 
every sheet contains exactly one nilpotent orbit.  
Sufficient conditions are given in \cite{BK} in order to ensure that the multiplicity with which an irreducible module occurs in the $G$-module decomposition of the ring of regular functions on an orbit is preserved along sheets. For instance, this always holds for ${\mathfrak sl}_n({\mathbb C})$. 

By a result of Katsylo \cite{katsylo}, a sheet $S$  can be described by means of the Slodowy slice through a nilpotent element $e$ lying in it and there exists a geometric quotient  $S/G$ for any sheet. 

The case of ${\mathfrak sl}_n({\mathbb C})$ is particularly clear: the sheets are disjoint (\cite{dix}), smooth, their $G$-orbits can be described by a quotient that is an affine space (\cite{kraft}), every sheet is a Dixmier sheet (\cite{OW}).

The interest in sheets in a Lie algebra has increased in recent years:  a class of nilpotent elements in a Lie algebra that are contained in a unique sheet and for which the sheet is smooth has been constructed in \cite{panyushev}, it has been proved that sheets in classical Lie algebras are always smooth in \cite{imhof} and 
the dimension of the sheets in a complex reductive Lie algebra has been  computed in \cite{moreau}.

The present paper addresses the analysis of sheets for a connected reductive algebraic group $G$ over an algebraically closed field $k$ of good characteristic acting on itself by conjugation. 
We will show affinities and discrepancies with some results holding for the adjoint action. 
The main goal of the paper is to show that sheets of conjugacy classes  are in bijection with $G$-conjugacy classes of triples $(M, tZ^\circ, \OO)$, with $M$ the connected centralizer of a semisimple element, $Z$ its center, $tZ^\circ$ a coset in $Z/Z^\circ$ with the property  that $M$ is the centralizer of $tZ^\circ$, and $\OO$ is a unipotent $M$-conjugacy class that is itself a sheet in $[M,\,M]$ (Theorem \ref{triplets} (1)). Another form of this result is that sheets of conjugacy classes are in bijection with pairs $(L,\,\OO)$ with $L$ a Levi factor of a parabolic subgroup of $G$ (from now on such an $L$ will be called a Levi subgroup) and $\OO$ a conjugacy class in $[L,\,L]$
which is itself a sheet (rigid conjugacy class). In this form the similarity with the Lie algebra case is evident.

When $G=PGL_n({\mathbb C})$ the connected centralizer $M$ of a semisimple element is always a Levi subgroup and $Z(M)$ is connected, so the parametrization coincides with the parametrization of sheets in ${\mathfrak sl}_n$.  However, for $G$ with components that are not of type $A_n$, we have, together with sheets that have a counterpart in the corresponding Lie algebra, also sheets that have a completely different behaviour. Since, for $G$ simple,  conjugacy classes of connected centralizers are described in \cite{sommers}, cosets in their center modulo the connected center have been analyzed in \cite{mcninch-sommers}, and rigid nilpotent orbits are listed in \cite{ke,el, sp}, our result allows a complete classification of sheets of conjugacy classes, which is part of a work in progress.

Two new phenomena occur which are two faces of the same token: the connected centralizer of a semisimple element in $G$ is not always a Levi subgroup, and there are noncentral, non-unipotent elements whose $G$-orbit is a whole sheet. The key idea in order to deal with the connected centralizer of a semisimple element that is not a Levi subgroup is to replace it by the minimal Levi subgroup containing it. Then one may still use induction of orbits as in \cite{bo, lusp}. It turns out that non-unipotent induced conjugacy classes may be described in terms of an induced unipotent conjugacy class in the centralizer of the semisimple part of a representative (Proposition \ref{crucial}). 

We also prove that a semisimple element lies, as for the case of Lie algebras, in a unique sheet (Theorem \ref{triplets} (4)).  Sheets containing a semisimple element correspond to triples with $\OO=\{1\}$.
Contrarily to what happens in the Lie algebra case, where a sheet always contains a nilpotent orbit,  
a sheet in the group contains a unipotent conjugacy class (up to a central element) if and only if the subgroup $M$ in its corresponding triple is a Levi subgroup.  Moreover, it contains a genuine unipotent element if and only if $tZ^\circ=Z^\circ$ (Theorem \ref{triplets} (2)). In those cases, the unipotent conjugacy class involved is unique (Theorem \ref{triplets} (3)).  So, when $M$ is not a Levi subgroup we cannot expect to have a straightforward analogue of Katsylo's result. A generalization of this result as a parametrization of conjugacy classes in a sheet in terms of a geometric quotient is part of a forthcoming project.

Several properties and invariants of orbits are preserved along sheets. For instance, in the general setup it has been shown in \cite{arzha} that complexity, i.e., the minimal codimension of orbits for a Borel subgroup $B$ in a $G$-homogeneous space is constant along sheets. For this reason, we view the understanding of sheets as part of a program in the comprehension of conjugacy classes.
We expect that the analysis of sheets will have applications to the study of the intersections of  conjugacy classes with Bruhat cells, to the analysis of the combinatorics of closures of orbits for the action of a Borel subgroup on conjugacy classes, and, in special cases, to the $G$-module decomposition of the ring of regular functions on a conjugacy class.

It is worthwhile to notice that $G$-conjugacy classes of triples $(M, tZ^\circ, M\cdot u)$ with $M$ and $tZ^\circ$ as in our case and 
 $M\cdot u$ a distinguished unipotent conjugacy class in $M$, have been used in \cite{mcninch-sommers} in order to describe conjugacy classes in the component group $A(u)$ of the centralizer of a unipotent element $u$ in good characteristic.

\section{Notation}

Unless otherwise stated, $G$ is a connected, reductive algebraic group over an algebraically closed field $k$ of good characteristic, i.e., not bad for any simple component of $[G,G]$. For the definition of good primes, see \cite[\S I.4.3]{131}.

Let $T$ be a fixed maximal torus of $G$ and let $\Phi$ be the associated
root system. Let $B\supset T$ be a Borel subgroup with unipotent radical $U$, let
$\Delta=\{\alpha_1,\ldots,\alpha_n\}$ be the basis of $\Phi$ relative
to $(T,\,B)$.  If $\Phi$ is irreducible, we shall denote by $-\alpha_0$ the highest positive root in $\Phi$. In this case we write $\alpha_0=\sum_{i=1}^n c_i\alpha_i$ and we set $c_0:=1$. 
The numbering of simple roots is as in \cite{bourbaki}.
The Weyl group of $G$ will be denoted by $W$.
Let $V$ be a variety and let $x\in V$; we shall denote by $V_x$ the connected component of $V$ containing $x$ so that, if $V$ is an algebraic group,  we have $V_1=V^\circ$. 
The centralizer of an element $x\in G$ in a subgroup $H$ of $G$ will be denoted by $H^x$ and its identity component will be denoted by $H^{x\circ}$. 
Let $G$ act regularly on an irreducible
variety $X$.
For $n\geq 0$, we shall denote by $X_{(n)}$ the locally closed subset  $X_{(n)}=\{x\in X~|~ \dim G\cdot x=n\}$. 
For a subset $Y\subset X$, if $m$ is the maximum integer $n$ for which $Y\cap X_{(n)}\neq\emptyset$,  the open subset $X_{(m)}\cap Y$ will be denoted by $Y^{reg}$.  We will investigate the case in which $X=G$ and the action is by conjugation. 
If $s\in H^{reg}$ for a subgroup $H$ of $G$, then $(H^{reg})_s$ is well defined, we have $(H^{reg})_s=(H_s)^{reg}$ and we will denote it by $H^{reg}_s$.

\section{Pseudo-Levi subgroups and Levi envelopes}

A pseudo-Levi subgroup  of $G$ is a subgroup of the form $G^{s\circ}$ for some semisimple element $s\in G$.  Several results on pseudo-Levi subgroups are to be found in \cite{sommers, mcninch-sommers}. By \cite[Exercise I.4.7]{131}, if the characteristic  of the base field is good for $G$ then it is good for any of its pseudo-Levi subgroups.

\begin{remark}\label{deriziotis}(Deriziotis' criterion) Let  us assume that $\Phi$ is an irreducible root system and let $\Psi\subset\Phi$ be a root subsystem. The   
subgroup $H=\langle T,\,X_\alpha,\,\alpha\in\Psi\rangle$ is a pseudo-Levi subgroup if and only if 
$\Psi$ has a basis that is $W$-conjugate to a subset of $\Delta\cup\{\alpha_0\}$, where $-\alpha_0$ is the highest root in $\Phi$. The original proof of Deriziotis holds when the base field is the algebraic closure of a finite field. A general proof  in good characteristic is in \cite[Proposition 30, Proposition 32]{mcninch-sommers}, see also \cite{sommers},\cite[\S 5.5]{lusztig}. 
\end{remark}

 \begin{example}\label{exc}Let $G=Sp_4(k)$ be the group of $4\times 4$ matrices  that leave invariant the bilinear form whose matrix with respect to the canonical basis in $k^4$ is $\left(\begin{array}{cc}
0&-I_2\\
I_2&0
\end{array}\right)$. Then the centralizer of the diagonal matrix $s={\rm diag}(-1,1,-1,1)$ is of type $A_1\times A_1$, hence it cannot be a Levi subgroup. It corresponds to the subset $\{\alpha_0,\,\alpha_2\}$ of $\{\alpha_0\}\cup\Delta$.
 \end{example} 

\begin{remark}\label{slodowy}{(Slodowy's criterion) Let $\Phi$ be the root system of $G$ with respect to a fixed maximal torus $T$. We say that a subset $\Sigma\subset  \Phi$ is ${\mathbb Q}$-closed if ${\rm span}_{\mathbb Q}(\Sigma)\cap \Phi=\Sigma$. According to \cite[Corollary 3.5]{slodowy},  if $s\in T$ is a semisimple element and $\Sigma=\{\alpha\in \Phi~|~\alpha(s)=1\}$,  
the pseudo-Levi subgroup $G^{s\circ}=\langle T,\,X_\alpha,\,\alpha\in\Sigma\rangle$ is a Levi subgroup if and only if $\Sigma$ is ${\mathbb Q}$-closed. }
\end{remark}

%

In the sequel we will often make use of the following sets, for $s\in T$:

\begin{equation}\label{zeta}Z(G^{s\circ})=\{g\in G~|~G^{g}\supset G^{s\circ}\}=\{g\in T~|~G^{g}\supset G^{s\circ}\}=\{g\in T~|~G^{g\circ}\supset G^{s\circ}\}\end{equation}
\begin{equation}\label{zetareg}Z(G^{s\circ})^{reg}=\{g\in T~|~G^{g\circ}\supset G^{s\circ}\}^{reg}=\{g\in T~|~G^{g\circ}=G^{s\circ}\}.\end{equation}

%
%

\bigskip

\begin{remark}\label{immagini}
Let $\pi\colon G\to \overline{G}$ be a central isogeny of connected and reductive groups. It is not hard to verify that, if $s\in T$ and $\pi(s)=\overline{s}$ then
$\pi(G^{s\circ})=\overline{G}^{\overline{s}\circ}$. Moreover, the descriptions in \eqref{zeta}, \eqref{zetareg} show that $\pi(Z(G^{s\circ}))=Z(\overline{G}^{\overline{s}\circ})$ and 
$\pi(Z(G^{s\circ})^{reg})=Z(\overline{G}^{\overline{s}\circ})^{reg}$. Since $\pi(Z(G^{s\circ})^\circ)$ is a closed, connected subgroup of $Z(\overline{G}^{\overline{s}\circ})$ of finite index, we also have
$\pi(Z(G^{s\circ})^\circ)=Z(\overline{G}^{\overline{s}\circ})^\circ$.
\end{remark}

The next Proposition is a reformulation of some results in \cite{mcninch-sommers}.

\begin{proposition}\label{simply}Let $G$ be connected and reductive, let $M$ be a pseudo-Levi subgroup 
and let $Z=Z(M)$. The following are equivalent:
\begin{enumerate}
\item $M$ is a Levi subgroup.
\item $Z=Z(G)Z^\circ$. 
\item $Z^{reg}=Z(G)(Z^\circ\cap Z^{reg})$. 
\item $Z^\circ\cap Z^{reg}\neq\emptyset$.
\item $(Z^\circ)^{reg}=Z^{reg}\cap Z^\circ$.
\item For every irreducible component $Z_r$ of $Z$ we have 
$Z_r\cap Z^{reg}\neq\emptyset$.
\end{enumerate}
\end{proposition}
\pf Let $s$ be such that $M=G^{s\circ}$. It is not restrictive to assume that $s\in T$.

We first show that 1. implies 2. under the assumption that $G$ is simple and adjoint.
By Remark \ref{deriziotis} we may assume that $M=\langle T,\,X_\alpha,\,\alpha\in\Psi\rangle$ for $\Psi$ generated by a subset $J$ of $\tilde{\Delta}=\Delta\cup\{\alpha_0\}$.
By \cite[Lemma 33]{mcninch-sommers} or \cite[\S 2.1]{sommers}  the group $Z/Z^\circ$ is cyclic of order $d=gcd(c_i~|~\alpha_i\in \tilde{\Delta}\setminus J)$. In particular, if 1. holds then $J\subset\Delta$ and since $c_0=1$ we have 2. Let now $G$ be simple and let  $\pi\colon G\to \overline{G}$ be the projection onto the adjoint
group and let $\overline{s}=\pi(s)$. If $M$ is a Levi subgroup then
$\pi(M)=\overline{G}^{\overline{s}\circ}$ is a Levi subgroup
of $\overline{G}$.
By Remark \ref{immagini} and the adjoint case we have $\pi(Z)=Z(\overline{G}^{\overline{s}\circ})=Z(\overline{G}^{\overline{s}\circ})^\circ=\pi(Z^\circ)$, hence $Z={\rm Ker}(\pi)Z^\circ=Z(G)Z^\circ$.
For the general case we consider the central isogeny $Z(G)^\circ\times G(1)\times\cdots\times G(m)\to G$ where $G(1),\ldots,\,G(m)$ are the components of $G$ and we use Remark \ref{immagini} once more.

If 2. holds then 3. follows by taking the intersection with
$Z^{reg}$.  If 3. holds then 4. is clearly true. 
If  $t\in Z^\circ \cap Z^{reg}$, then  $M=G^{t\circ}=C_G(tZ^\circ)=C_G(Z^\circ)$  hence $M$ is the connected centralizer of a torus. By \cite[\S II.4.1]{131}, $M$ is generated by the maximal torus containing it and the root subgroups corresponding to the roots vanishing on it. By
\cite[Proposition III.8.9]{borel}, the corresponding root subsystem is ${\mathbb Q}$-closed hence  4. implies 1.  

Since $\max\dim_{x\in Z}G\cdot
x\geq \max\dim_{x\in Z^\circ}G\cdot x$, we have equivalence of
4. and 5.

If 2 holds then 4 and 5 hold and for every component $Z_r$ we
have
$Z_r=zZ^\circ$ for some $z\in Z(G)$, and
$(zZ^\circ)^{reg}=z(Z^\circ)^{reg}$ so
6. follows. Conversely, 6. with $r=1$ is 4.
\hfill$\Box$

\bigskip

\begin{definition}Given a pseudo-Levi $H$ with center $Z$, the Levi-envelope $L$ of $H$ in $G$ is $L=C_G(Z^\circ)$. 
\end{definition}

\begin{lemma}\label{zetagiesse}Let $H$ be a pseudo-Levi subgroup containing $T$ and let $L$ be its Levi envelope in $G$. Then $L$ is the minimum Levi subgroup of $G$ containing $H$ and $Z(H)^\circ=Z(L)^\circ$. Moreover, 
if $\Sigma$ is the root subsystem of $H$ relative to $T$, then  $\overline{\Sigma}={\rm span}_{\mathbb Q}(\Sigma)\cap\Phi$, the ${\mathbb Q}$-closure of 
$\Sigma$ in $\Phi$ is the root subsystem of $L$ relative to $T$.
\end{lemma}
\pf Let $s\in T$ such that $G^{s\circ}=H$. Let $L'$ be a Levi subgroup of $G$ containing $H$.  It is not hard to verify that $L'=C_G(T')$ for some torus $T'$. Then $T'\subset Z(H)^\circ$ because $s\in L$ so $T'\subset G^{s\circ}$ and centralizes $G^{s\circ}$. 
Thus, $L\subset L'$.  Besides, $Z(H)^\circ\subset Z(L)^\circ$ by construction. On the other hand, $H\subset L$ so $Z(L)^\circ\subset Z(H)^\circ$. 

The Levi subgroups of $G$ containing $H$ contain $T$, so they correspond to ${\mathbb Q}$-closed root subsystems of $\Phi$ containing $\Sigma$ 
and  $\overline{\Sigma}$ is clearly the minimum such root subsystem.\hfill$\Box$  

\begin{example}Let $G$ and $s$ be as in Example \ref{exc}. The Levi envelope $L$ of  $G^{s\circ}$ is $G$. Indeed, $\alpha_1=\frac{1}{2}(\alpha_2-\alpha_0)$ so the root subgroups $X_{\pm\alpha_1}$ and $X_{\pm\alpha_2}$  lie in $L$.
\end{example}

The following Lemmas will be needed in the sequel.

\begin{lemma}\label{pseudo-levi}A Levi subgroup $L$ of a pseudo-Levi subgroup $M$ of $G$ is a pseudo-Levi subgroup of $G$.
\end{lemma}
\pf The statement easily follows once we have proven it for $G$  simple. We may assume that $T\subset L\subset M$. By Remark \ref{deriziotis}  the root system of $M$ relative to $T$ has a basis that is $W$-conjugate to a subset $J$ of $\{\alpha_0\}\cup\Delta$ and the root system of $L$ relative to $T$ has a basis that  is $W$-conjugate to a subset of $J$. Applying Remark \ref{deriziotis} once more we have the statement.\hfill$\Box$

\begin{remark}In general a pseudo-Levi subgroup $M$ of a pseudo-Levi subgroup $N$ of $G$ is not pseudo-Levi in $G$. For instance, if $G={\rm Sp}_8(k)$ is the group of $8\times 8$ matrices  that leave invariant the bilinear form whose matrix with respect to the canonical basis in $k^8$ is $\left(\begin{array}{cc}
0&-I_4\\
I_4&0
\end{array}\right)$, then if $M$ is the centralizer of $s={\rm diag}(-I_2,I_2,-I_2,I_2)$, the centralizer $N=M^{r\circ}$  of $r={\rm diag}(-1,1,-1,1,-1,1-1,1)$ in $M$, is of type $A_1\times A_1\times A_1\times A_1$, so it cannot be a pseudo-Levi subgroup of $G$ by Remark \ref{deriziotis}.
\end{remark}

\begin{lemma}\label{pinco}Let $N$ be a Levi subgroup of a pseudo-Levi subgroup $M$ of $G$ and let $L$ be be its Levi-envelope in $G$. Then $L\cap M$ is  the Levi-envelope of $N$ in $M$.
\end{lemma}
\pf Let $Z=Z(N)$. We have $L=C_G(Z^\circ)$. Then $L\cap M=C_G(Z^\circ)\cap M=C_M(Z^\circ)$ and we have the statement.
\hfill$\Box$

\section{Jordan classes}

In this section we will introduce the notion of Jordan classes, which, just as in the case of the adjoint action of a group on its Lie algebra, is crucial for our purposes.
In analogy to \cite{bo,BK}, we define the following equivalence relation on a reductive group $G$:
for $g\in G$ with Jordan decomposition $g=su$ we have $g\sim h$ if there exists $x\in G$ such that $x^{-1}hx=rv$ with $G^{s\circ}=G^{r\circ}$, $r\in Z(G^{s\circ})_s$ and $G^{s\circ}\cdot u=G^{s\circ}\cdot v$.

\begin{example}Let $SO_{8}(k)$ be the group of $8\times 8$ matrices  of determinant equal to$1$ that leave invariant the bilinear form whose matrix with respect to the canonical basis in $k^8$ is $\left(\begin{array}{cc}
0&I_4\\
I_4&0
\end{array}\right)$. The two matrices  $s={\rm diag}(-1,I_3,-1,I_3)$ and $r={\rm diag}(t,I_3,t^{-1},I_3)$ for some $t\neq0,\,\pm1$ are equivalent. 
\end{example}

\begin{proposition}\label{equivalence}
Let $g\in G$ have  Jordan decomposition $g=su$. Then, for $h\in G$ we have $g\sim h$ if and only if  $G\cdot h\cap
Z(G^{s\circ})^{reg}_su\neq\emptyset$.
\end{proposition}
\pf If $g\sim h$ then there exixts $rv\in G\cdot h$  such that $G^{r\circ}=G^{s\circ}$ with $r\in Z(G^{r\circ})_r=Z(G^{s\circ})_s$ and $u=x\cdot v$ for some $x\in G^{s\circ}$. Then, $x\cdot rv=ru\in G\cdot h \cap  Z(G^{s\circ})^{reg}_su$.

Conversely, if $ru\in G\cdot h\cap Z(G^{s\circ})^{reg}_su$ then $r\in Z(G^{s\circ})_s$ with $\dim G^s=\dim G^r$. So, $G^{s\circ}\subset G^r$ and by the dimension condition we conclude that $G^{s\circ}=G^{r\circ}$. Thus, $h\sim g$.\hfill$\Box$

\bigskip

The equivalence classes with respect to  $\sim$  are
called Jordan classes. The Jordan class of $g=su$ in $G$ will be denoted by $J_G(g)$ and by  Proposition \ref{equivalence} we have 
$J_G(g)=G\cdot (Z(G^{s\circ})^{reg}_su)$.  Jordan classes are irreducible.
There are only finitely many of them in a group $G$. Indeed, one may always assume that $s\in T$ 
so that $G^{s\circ}$ is determined by a root
subsystem of $\Phi$, and there are thus only finitely many possible $G^{s\circ}$. Moreover,  each $Z(G^{s\circ})$ has only finitely many irreducible
components and each $G^{s\circ}$ has only finitely many unipotent conjugacy
classes. 

Is it not hard to verify that if $g=su\in G_{(n)}$ then $J_G(g)\subset G_{(n)}$ and $\overline{J_G(g)}\subset  \bigcup_{m\leq n}G_{(m)}$  since $\overline{G_{(n)}}=\bigcup_{m\leq n}G_{(m)}$. 

Our first  goal  is to understand the closure of a Jordan class.
 A fundamental notion in the comprehension of closures of Jordan classes for the adjoint action of a group on its Lie algebra is played by induction of orbits (\cite{bo, lusp}). 
This notion can be adapted to our situation. 
 
Let  $M$ be a Levi subgroup in a connected reductive group $K$, let $P=MU_P$ be a parabolic subgroup of $K$ with unipotent radical $U_P$ and let $\OO_1$ be a conjugacy class in $M$.
We define the {\em conjugacy class induced by $\OO_1$} as $\OO:={\rm Ind}_{M,P}^K(\OO_1):=K\cdot(\OO_1 U_P)^{reg}$. One proves as in \cite[\S 2.1]{bo} that $\OO$ is indeed a $K$-conjugacy class. 
Induced unipotent conjugacy classes have been extensively studied in \cite{lusp}.  It is well-known that they are independent  of the choice of the parabolic subgroup $P$.

\begin{lemma}\label{chiusura-indotte}Let $P=LU_P$ be a parabolic subgroup of a connected reductive group $G$ and let $\OO$ be an $L$-conjugacy class. Then $\overline{{\rm Ind}_{L,P}^G(\OO)}=G\cdot (\overline{\OO}U_P)$.
\end{lemma}    
 \pf It can be proved as in \cite[Lemma p. 290]{bo}.
\hfill$\Box$

\bigskip

\begin{lemma}\label{semis-part1}Let $P=LU_P$ be a parabolic subgroup of a connected reductive group $G$ and let $su\in L$.
Then the semisimple parts of elements in $su U_P$ are $P$-conjugate to $s$.\end{lemma}    
 \pf This is proved as in \cite[\S 2.1]{bo}. 
\hfill$\Box$

\begin{lemma}\label{representative}Let $P=LU_P$ be the Levi decomposition of a parabolic subgroup of $G$. Let $g=su$ be the Jordan decomposition of an element $g\in L$.
Then, there exists a representative of ${\rm Ind}_{L,P}^G(L\cdot su)$ lying in $suU_P^s$.\end{lemma}
\pf The closed set $suU_P$ is $U_P$-stable and it contains $Y=U_P\cdot(su U_P^s)$. Moreover, if $\pi\colon U_P\times su U_P^s\to Y$ is the natural morphism, then for every $y$ in an open subset $V$ of $Y$  we have $\dim Y=\dim U_P+\dim U_P^s-\dim \pi^{-1}(y)$. We can make sure that $V\cap suU_P^s\neq\emptyset$ and that $y$ lies in this intersection. 
Then,
$$\pi^{-1}(y)=\{(v,x)\in U_P\times su U_P^s~|~v\cdot x=y\}.$$
Every $x$ in $su U_P^s$ has semisimple part equal to $s$, so this holds for $y$, too. Therefore, any $v\in U_P$ for which $v\cdot x=y$ must preserve $s$, hence $U_P\cdot y\cap suU_P^s=U_P^s\cdot y$ and
$$\pi^{-1}(y)=\{(v,x)\in U_P^s\times U_P^s\cdot y~|~v\cdot x=y\}=\{(v,\,v^{-1}\cdot y)\in U_P^s\times U_P^s\cdot y\}\cong U_P^s.$$
Hence, $\dim Y=\dim U_P$ and $\overline{Y}=suU_P$. Thus, $Y\cap(suU_P)^{reg}\neq\emptyset$.\hfill$\Box$

\bigskip

\begin{proposition}\label{crucial}Let $L$ be a Levi subgroup of $G$ containing $T$ and let $P$ be any parabolic subgroup of $G$ with Levi decompostion $P=LU_P$. If $L\cdot su$ is an $L$-conjugacy class with $s\in T$ and $u\in L^{s\circ}$ we have
\begin{equation}\label{indotta}{\rm Ind}_{L,P}^G(L\cdot su)=G\cdot(s\,{\rm Ind}_{L^{s\circ}}^{G^{s\circ}}(L^{s\circ}\cdot u)).\end{equation}
The class ${\rm Ind}_{L,P}^G(L\cdot su)$ is independent of the choice of $P$, it depends only on $L^{s\circ}$ and not on $L$,  the semisimple part of an element in it is conjugate to $s$ and 
\begin{equation}\label{dimensione}\dim {\rm Ind}_{L,P}^G(L\cdot su)=\dim G-\dim L+\dim L\cdot su.\end{equation} In particular, if $L\supset G^{s\circ}$ then ${\rm Ind}_{L,P}^G(L\cdot su)=G\cdot su$.
\end{proposition}
\pf Let $\Psi$ and $\Psi_s$ be the root systems of $L$ and $G^{s\circ}$, respectively,  relative to $T$.  Then the root system of the subgroup $L^{s\circ}$ is $\Psi\cap \Psi_s$ and it is ${\mathbb Q}$-closed in $\Psi_s$, hence $L^{s\circ}$ is a Levi subgroup of $G^{s\circ}$ by Remark \ref{slodowy}. Moreover, $T\subset P$ so there is a Borel subgroup $B$ of $G$ for which $T\subset B\subset P$. Then, $B^s=G^{s\circ}\cap B$ is a Borel subgroup of $G^{s\circ}$. Thus, $P^{s\circ}\supset B^{s}$ is a parabolic subgroup of $G^{s\circ}$. Since $L$ and $U_P$ are $T$-stable and $s\in T$, the Levi decomposition of $P$ induces the decomposition $P^{s\circ}=L^{s\circ}U_P^s$, hence $U_P^s$ is the unipotent radical of $P^{s\circ}$.  By Lemma \ref{representative} there exists $x\in {\rm Ind}_{L,P}^G(L\cdot su)\cap su U_P^s$. So, $(su U_P^s)^{reg}\subset (su U_P)^{reg}$ and 
$$s{\rm Ind}_{L^{s\circ}}^{G^{s\circ}}(L^{s\circ}\cdot u)=G^{s\circ}\cdot(su U_P^s)^{reg}\subset G\cdot (suU_P)^{reg}={\rm Ind}_{L,P}^G(L\cdot su),$$ whence formula \eqref{indotta}.
Since induced unipotent conjugacy classes are independent of the choice of the parabolic (\cite[Theorem 2.2]{lusp},\cite[Satz 2.6]{bo}), we deduce the same statement for general induced conjugacy classes. Besides, 
$$\begin{array}{rl}
\dim{\rm Ind}_{L,P}^G(L\cdot su)&=\dim G-\dim G^{s\circ}+\dim{\rm Ind}_{L^{s\circ}}^{G^{s\circ}}(L^{s\circ}\cdot u)\\
&=\dim G-\dim L^{s\circ}+\dim L^{s\circ}\cdot u\\
&=\dim G-\dim L+\dim L\cdot su
\end{array}$$
where we have used the formula for induced unipotent conjugacy classes in \cite[Theorem 1.3(a)]{lusp},\cite[Satz 3.3]{bo}.
The remaining statements easily follow from \eqref{indotta}.
\hfill$\Box$

\bigskip

As a consequence of Proposition \ref{crucial} we will write ${\rm Ind}_L^G$ instead of ${\rm Ind}_{L,P}^G$.

\begin{remark}\label{semis-part}If $G\cdot
  h\subset \overline{G\cdot g}$, the semisimple part of $h$ is conjugate to the
  semisimple part of $g$. Indeed, the subset of classes with semisimple part lying in a fixed conjugacy class is closed because it is a fiber of the map $G\to G/\!/G$ of $G$ to its categorical quotient.
\end{remark}

\begin{proposition}\label{chiusura}Let $G$ be a connected reductive group and let $J_G(g)$ be the Jordan class of the element $g$ with Jordan decomposition $g=su$, with $s\in T$. Let $Z=Z(G^{s\circ})$, let $L=C_G(Z^\circ)^\circ$ be the Levi envelope of $G^{s\circ}$ in $G$. Then 
\begin{equation}\label{closure}\overline{J_G(g)}=\bigcup_{z\in{Z^\circ}}\overline{{\rm Ind}_{L}^G(L\cdot zsu)}\end{equation}
and
\begin{equation}\label{reg-closure}\overline{J_G(g)}^{reg}=\bigcup_{z\in{Z^\circ}}{\rm Ind}_{L}^G(L\cdot zsu).\end{equation}
\end{proposition}
\pf The proof follows the line of \cite{bo, BK}. By Lemma \ref{zetagiesse} we have $Z^\circ=Z(L)^\circ$. Let   $P=LU_P$ be the Levi decomposition of  a parabolic subgroup with Levi component $L$. Let $R=Z^\circ\overline{L\cdot su}U_P$. Since $s,\,u\,\in G^{s\circ}\subset L$ we see that $R$ is $P$-stable. Moreover, $R$ is closed. Indeed, since $L\cap U_P=1$ it is enough to show that $ 
Z^\circ\overline{L\cdot su}$ is closed. As $L$ is reductive and connected and $[L,\,L]$ is closed,
$Z^\circ \overline{L\cdot su}=Z^\circ \overline{[L,\,L]\cdot s'u'}$ for some $s'u'\in[L,\,L]$ and the latter is closed because it is the image of a closed set under the isogeny $Z(L)^\circ\times [L,\,L]\to L$. By \cite[Chapter II.13, Lemma 2]{steinberg}, the saturation
$G\cdot R$ is closed. We have
$$\begin{array}{rl}
J_G(su)&=G\cdot (Z_s^{reg}u)\subset G\cdot(Z^\circ su)=G\cdot P\cdot(Z^\circ su)\subset G\cdot R
\end{array}$$
thus $\overline{J_G(su)}\subset \overline{G\cdot R}=G\cdot R$. 
 
Let us consider the set $M=P\cdot (Z^{reg}_{s}u)\subset J_G(g)$. 

For every $zu\in Z^{reg}_{s}u\subset L$ we have $U_P\cdot zu\subset zu U_P$. Moreover, $\dim U_P\cdot zu=\dim U_P-\dim U_P\cap G^{z}\cap G^u$ and
$\dim U_P\cap G^{z}\cap G^u=\dim (U_P\cap G^{s}\cap G^u)^\circ\leq \dim U_P\cap G^{s\circ}\leq \dim U_P\cap L=0$. So $U_P\cdot zu= zu U_P$. 
Then $M=L\cdot(Z^{reg}_{s}u)U_P\subset L\cdot(Z^\circ su)U_P\subset R$. Since 
$$\overline{Z^{reg}_{s}u}=\overline{Z^{reg}_{s}}\,u=Z_{s}u=Z^\circ su$$ we have 
$Z^\circ su U_P= \overline{Z^{reg}_{s}uU_P}\subset \overline{M}$. By 
$P$-stability of $\overline{M}$ we have 
$$L\cdot(Z^\circ su U_P)=Z^\circ (L\cdot su)U_P\subset\overline{M},$$ thus
$$\overline{Z^\circ (L\cdot su)U_P}=\overline{Z^\circ (L\cdot su)}U_P=R\subset\overline{M}$$ so $\overline{M}=R$. By $G$-stability 
$G\cdot R=G\cdot\overline{M}\subset \overline{J_G(g)}$, whence $\overline{J_g(g)}=G\cdot(Z^\circ\overline{L\cdot su}U_P)$. 
Formula \eqref{closure}  follows from Lemma \ref{chiusura-indotte}, as
$$G\cdot (Z^\circ\overline{L\cdot su}U_P)=\bigcup_{z\in Z^\circ}G\cdot (\overline{L\cdot zsu}U_P)=\bigcup_{z\in Z^\circ}\overline{{\rm Ind}_{L}^G(L\cdot zsu)}.$$
Equation \eqref{reg-closure} follows from  \eqref{dimensione} and the fact that $L\cdot zsu=z L\cdot su$ since $z\in Z^\circ$.
\hfill$\Box$

\bigskip

\begin{lemma}\label{levi}Let $g=su\in \overline{J(x)}^{reg}$. If $x$ has Jordan decomposition $x=rv$ with $s\in Z^\circ r$ for $Z=Z(G^{r\circ})$, then
$G^{r\circ}$ is a Levi subgroup of $G^{s\circ}$. (It is always possible to find a representative of the Jordan class of this form).
\end{lemma}
\pf Clearly $G^{r\circ}\subset G^{s\circ}$. 
%
Then, $G^{r\circ}=C_G(Z^\circ r)^\circ=C_G(Z^\circ s)^\circ=(G^{s\circ}\cap C_G(Z^\circ))^\circ=C_{G^{s\circ}}(Z^\circ)$. Since $Z^\circ$ is a torus in $G^{s\circ}$, we have the statement.
\hfill$\Box$

%
%

\bigskip

\begin{proposition}\label{closure-inclusion}The closure of a Jordan class  in $G$ is a union of  (closures of ) Jordan classes.
\end{proposition}
\pf Let $g\in G$ with Jordan decomposition $g=su\in TU$ and let $x\in \overline{J_G(g)}$. We shall show that $J_G(x)$ lies in $\overline{J_G(g)}$.  Since $G\cdot x\subset \overline{J_G(g)}$  we may replace $x$ by a suitable $G$-conjugate. With same notation as in Proposition \ref{chiusura}, we may assume that $x\in \overline{L\cdot zsu}U_P$ for some $z\in Z^\circ$.
By Remark \ref{semis-part} applied to $L\cdot zsu$, we may take $x\in zsvU_P\subset (L\cdot zsv)U_P\subset \overline{L\cdot zsu}U_P$ for some unipotent $v\in G^{zs\circ}$. 
Applying Lemma \ref{semis-part1},  we see that the semisimple part of $x$ is $P$-conjugate to $zs$. Thus, there is some $y\in P\cdot x\cap zsvU_P\subset \overline{L\cdot zsu}U_P$ with Jordan decomposition $y=(zs) u'$ with $u'\in vU_P$. Then,
$G^{zs\circ}\supset G^{s\circ}$ so $Z(G^{zs\circ})^\circ\subset Z^\circ$,  $Z(G^{zs\circ})_{zs}\subset Z_s$, and 
$$Z(G^{zs\circ})^{reg}_{zs}u'\subset Z(G^{zs\circ})^\circ y\subset Z^\circ  svU_P\subset Z^\circ \overline{L\cdot zsu}U_P.$$
Hence, $J_G(x)=J_G(y)\subset \overline{J_G(g)}$.\hfill$\Box$

\bigskip

As a consequence of Proposition \ref{closure-inclusion}, the set ${\mathcal J}$ of Jordan classes in $G$ has a natural  partial order:
$J_G(g)\leq J_G(h)$  if $J_G(g)\subset \overline{J_G(h)}^{reg}$. 

 \begin{remark}\label{generates}Let $G$ be a simple group of adjoint type. Let $M$ be a pseudo-Levi subgroup of $G$ containing $T$ and let $\Delta'\subset\{\alpha_0\}\cup\Delta$ be a basis of its root system. Let $Z=Z(M)$ and let $t\in Z$. 
With the same arguments as in \cite[Proposition 15(2)]{mcninch-sommers} one shows that
$tZ^\circ\cap Z^{reg}\neq\emptyset$ if and only if $M=C_G(tZ^\circ)^\circ$. 
 By \cite[Lemma 33]{mcninch-sommers} the group $Z/Z^\circ$ is cyclic of order ${\rm gcd}(c_i~|~\alpha_i\in \tilde{\Delta}\setminus \Delta')$
 and by \cite[Lemma 34]{mcninch-sommers} we have $M=C_G(tZ^\circ)^\circ$ if and only if $tZ^\circ$ generates $Z/Z^\circ$. 
\end{remark}

\begin{proposition}\label{bijection}There is a bijective correspondence between ${\mathcal J}$ and the set of $G$-conjugacy classes of triples $(M,\,Z^\circ t, \OO_1)$ where $M$ is a pseudo-Levi subgroup of $G$ with center $Z=Z(M)$, the coset $Z^\circ t\in Z/Z^\circ$ is such that $C_G(Z^\circ t)^\circ =M$, and $\OO_1$ is a unipotent conjugacy class in $M$. Moreover, if $J_1,\, J_2\in{\mathcal J}$
correspond respectively to the triples $(M_1,\,Z^\circ _1t_1,\,\OO_1)$ and $(M_2,\,Z^\circ_2t_2,\,\OO_2)$, then $J_1\leq J_2$ if and only if there exists $g\in G$ such that
$g\cdot M_2$ is a Levi subgroup of $M_1$, the coset $Z_1^\circ t_1\subset g\cdot Z_2^\circ t_2$, and $\OO_1={\rm Ind}_{g\cdot M_2}^{M_1}(g\cdot \OO_2)$. 
\end{proposition}
\pf Let $(M,\,Z^\circ t,\OO)$ be any such triple and let $v\in \OO$. We may attach to it the Jordan class $G\cdot( (Z(M)^\circ t)^{reg}v)=J_G(sv)$ for any $s\in (Z(M)^\circ t)^{reg}$. 
So, the assignment
 $(M,\,Z^\circ t, M\cdot v)\mapsto G\cdot (Z(M)^\circ tv)$ determines a surjective map from the set of the above triples to ${\mathcal J}$ which is clearly constant on $G$-orbits with respect to simultaneous conjugation. 
Let us assume that $(M',\, Z(M')^\circ t',\,M'\cdot v')$ is a triple for which $G\cdot ((Z(M)^\circ t)^{reg}v)=G\cdot ((Z(M')^\circ t')^{reg}v')$. Let $s\in (Z(M')^\circ t')^{reg}$. Then there exists $g\in G$ such that $g\cdot sv'\in (Z(M)^\circ t)^{reg}v$. Thus, $g\cdot M'=M$, so $g\cdot (Z(M)^\circ t)=Z(M')^\circ t'$ and $g\cdot v'=v$.  Therefore, the map induced on $G$-orbits of such triples is a bijection. 

Let $u_1\in \OO_1$ and $u_2\in\OO_2$. By \eqref{reg-closure} and \eqref{indotta} we have $\overline{J_2}^{reg}=\bigcup_{z\in Z_2^\circ t_2}{\rm Ind}_L^G(L\cdot zu)$ where $L$ is the Levi-envelope of $M_2$ in $G$ and $u\in\OO_2$. If $J_1\leq J_2$, then conjugating by  $g\in G$ we may change representative of $J_2$  so that 
$g\cdot M_2\subset M_1$ and $Z(M_1)^\circ t_1\subset g\cdot(Z(M_1)^\circ t_2)$. Moreover, by Lemma \ref{levi}, the subgroup $g\cdot M_2$ is a Levi subgroup of $M_1$ so by Lemma \ref{pinco}, it coincides with $g\cdot L\cap M_1$. 

By \eqref{indotta} we see that  
$$J_1\subset \overline{J_2}^{reg}=\bigcup_{g\cdot z\in g\cdot(Z_2^\circ t_2)}{\rm Ind}_{g\cdot L}^G(g\cdot (L\cdot zu))=\bigcup_{t\in g\cdot(Z^\circ_2t_2)}G\cdot t{\rm Ind}_{g\cdot M_2}^{M_1}(g\cdot \OO_2)$$ 
and for our choice of representatives we have $$(Z(M_1)^\circ t_1)^{reg}u_1\subset \bigcup_{t\in g\cdot(Z^\circ_2t_2)}t{\rm Ind}_{g\cdot M_2}^{M_1}(g\cdot \OO_2)$$ so $\OO_1={\rm Ind}_{g\cdot M_2}^{M_1}(g\cdot \OO_2)$.

Conversely, assuming without loss of generality that $g=1$, i.e., that $M_2$ is a Levi subgroup of $M_1$, that $Z_1^\circ t_1\subset Z_2^\circ t_2$ and $\OO_1={\rm Ind}_{M_2}^{M_1}(\OO_2)$, if $s\in(Z_1^\circ t_1)^{reg}$ and $u_1\in \OO_1$ then 
$$G\cdot su_1=G\cdot s{\rm Ind}_{M_2}^{M_1}\OO_2\subset \overline{J_2}^{reg}$$ by \eqref{indotta}, \eqref{reg-closure} and Proposition \ref{closure-inclusion}. 
\hfill$\Box$

\section{Sheets}

A sheet for the $G$-action on itself by conjugation is an irreducible component of
 $G_{(n)}$.  A sheet is clearly $G$-stable.  
 

\bigskip

\begin{proposition}\label{quando}The map $J_G(g)\mapsto \overline{J_G(g)}^{reg}$ induces a bijection between maximal elements of ${\mathcal J}$ and sheets of $G$. A sheet is a union of Jordan classes. 
\end{proposition}
\pf The varieties $G_{(n)}$ are union of Jordan classes, so they are  finite unions of the irreducible closed sets $\overline{J}^{reg}$ for $J\in {\mathcal J}$ with $J\subset G_{(n)}$. Thus, the sheets in $G_{(n)}$ are precisely the closed sets of the form  $\overline{J}^{reg}$ with $J\subset G_{(n)}$ maximal with respect to the partial order in ${\mathcal J}$. The last statement follows from Proposition \ref{closure-inclusion}.\hfill$\Box$

%
%
%
%

\begin{remark}The closure of a sheet in a complex reductive Lie algebra is not always a union of sheets. Some counterexamples are to be found in (\cite{bulois-counter}). It will follow from the classification result in Theorem \ref{triplets} that similar counterexamples may be constructed for sheets in a reductive group.  
\end{remark}

%

\begin{proposition}\label{equivalenze}Let $g\in G$ with Jordan decomposition $g=su$. The following are equivalent:
\begin{enumerate}
\item $J_G(g)$ is maximal with respect to the partial ordering in ${\mathcal J}$.
\item $G\cdot su$ is not induced from any Levi subgroup of $G$ not containing $G^{s\circ}$.
\item $G^{s\circ}\cdot u$ is not induced from any proper Levi subgroup of $G^{s\circ}$.
\end{enumerate}
\end{proposition}
\pf If $G\cdot g={\rm Ind}_L^G(L\cdot sv)$ with $L$ not containing $G^{s\circ}$ then we may consider the Levi subgroup $M=L^{s\circ}$ of $G^{s\circ}$ and its Levi envelope $L'$. It follows from Lemma \ref{pinco} and Proposition \ref{crucial} that $G\cdot g={\rm Ind}_{L'}^G(L'\cdot sv)$. Besides, $s\in Z(M)$ and $Z(G^{s\circ})^\circ s\subset Z(M)^\circ s$. We have:
$$C_G(Z(M)^\circ s)^\circ=(C_G(Z(M)^\circ)\cap G^{s\circ})^\circ=(L'\cap G^{s\circ})^\circ=M$$ so by the arguments in \cite[Proposition 15 (2)]{mcninch-sommers} there exists $t\in Z(M)^\circ s$ such that $G^{t\circ}=M$. Thus, by \eqref{reg-closure},  we have the inclusion
$J_G(g)\subset\bigcup_{z\in Z(M)^\circ}{\rm Ind}_{L'}^G(L'\cdot ztu)=\overline{J_G(tv)}^{reg}$. Since $L$ does not contain $G^{s\circ}$, we have a proper inclusion $G^{t\circ}\subset G^{s\circ}$ so $J_G(g)\neq J_G(tv)$. Therefore, $J_G(g)$ is not maximal. Hence 1 implies 2. 

If   $G^{s\circ}\cdot u={\rm Ind}_M^{G^{s\circ}}(M\cdot v)$ is induced from a proper Levi subgroup $M$ of $G^{s\circ}$ it follows from Lemma \ref{pseudo-levi} that $M$ is a pseudo-Levi subgroup of $G$ and, for $L$ the Levi-envelope $L$ of $M$ we have $G\cdot su={\rm Ind}_L^G(L\cdot su)$ by Lemma \ref{pinco} and \eqref{indotta}. Since $M=L\cap G^{s\circ}$ is a proper subgroup of $G^{s\circ}$, the subgroup $L$ does not contain  $G^{s\circ}$ so 2 implies 3. 

If $G^{s\circ}\cdot u$ is not induced, $J_G(g)$ is maximal in ${\mathcal J}$ by Proposition \ref{bijection}.
\hfill$\Box$ 

It is shown in \cite{bo} that the nilpotent orbits in a semisimple Lie algebra that are not induced are precisely the rigid orbits, that is, those orbits that constitute by themselves a sheet. 


\begin{definition} A conjugacy class in a semisimple group $H$ is rigid if it is a sheet. 
\end{definition}

In the Lie algebra case, rigid orbits are necessarily nilpotent. In the case of a group it turns out that not every rigid conjugacy class is unipotent, but all rigid conjugacy class are exceptional in the following sense. 

\begin{definition}A semisimple element $s$ in $H$ is exceptional if $Z(H^{s\circ})^\circ=1$. An element in $H$ is called exceptional if its semisimple part is so. A conjugacy class is exceptional if its elements are exceptional. 
\end{definition}

This notion was introduced in  \cite[\S 7]{DCK} where exceptional elements played the role of nilpotent elements in the generalization of a result from the Lie algebra case to the group case.

We may state the classification result on sheets:

\begin{theorem}\label{triplets}Let $G$ be a connected and reductive group $G$ over an algebraically closed field of good characteristic. Then:
\begin{enumerate}
\item[(1)] The sheets of $G$ are in one-to-one correspondence with $G$-orbits of triples $(M, Z^\circ t, \OO)$
where $M$ is a pseudo-Levi subgroup of $G$ with centre $Z=Z(M)$, $Z^\circ t$ is a class in $Z/Z^\circ t$ such that $C_G(Z^\circ t)^\circ=M$, and $\OO$ is a rigid unipotent conjugacy class in $M$.
\item[(2)] The sheets containing a unipotent conjugacy class up to a central factor are precisely those for which $M$ is a Levi subgroup. The sheets containing a genuine unipotent conjugacy class are those for which $M$ is a Levi subgroup  and $Z^\circ t=Z^\circ$.
\item[(3)] If a sheet contains a unipotent conjugacy class up to a central element, the involved unipotent class is unique.
\item[(4)]  The sheets containing a semisimple element are those for which $\OO=\{1\}$ and every semisimple element lies in a unique sheet.
\item[(5)] The sheet associated with $(M,\,Z^\circ t,\,M\cdot u)$ lies in $G_{(n)}$ with $n=\dim G-\dim M+\dim M\cdot u$.
Its dimension is equal to $\dim G+\dim Z^\circ-\dim M^{u}$. 
\end{enumerate}
\end{theorem}
\pf (1) The first statement follows from Proposition \ref{bijection}, Proposition \ref{quando} and Proposition \ref{equivalenze}.


\noindent (2) If $zv\in \overline{J_G(sv)}^{reg}$ with $z\in Z(G)$, then $G^{s\circ}$ is a Levi subgroup of $G$ by Lemma \ref{levi}. 
Conversely, if $M=G^{s\circ}=L$ is a Levi subgroup of $G$, by Proposition \ref{simply} we have that $Z^\circ s=z Z^\circ$ for some $z\in Z(G)$. Then, by \eqref{reg-closure}
\begin{equation}\label{unipo}
\overline {J_G(sv)}^{reg}=\bigcup_{t\in Z(M)^\circ}{\rm Ind}_{M}^G(M\cdot ztv)
\supset z  {\rm Ind}_{L}^G(M\cdot u).
\end{equation}

Clearly, $\overline{J_G(sv)}^{reg}$ contains a genuine unipotent conjugacy class if and only if $z\in Z^\circ$.

\noindent (3) Is immediate from \eqref{unipo}. 

\noindent (4) Let $S$ be a sheet containing $s$ and let $S=\overline{J_G(rv)}^{reg}$. Then
$G\cdot s={\rm Ind}_L^G(L\cdot zrv)$, with $r\in T$, $z\in Z(G^{r\circ})^\circ$, $v\in G^{r\circ}$ and $L$  the Levi-envelope of $G^{r\circ}$ in $G$.
By Lemma \ref{levi}, the subgroup $G^{r\circ}$ is a Levi subgroup of $G^{s\circ}$ and we have $(L\cap G^{s\circ})=G^{r\circ}$ by Lemma \ref{pinco}. 
It follows from Lemma \ref{crucial} that $G\cdot s=G\cdot s{\rm Ind}_{G^{r\circ}}^{G^{s\circ}}(G^{r\circ}\cdot v)$. So
$1={\rm Ind}_{G^{r\circ}}^{G^{s\circ}}(G^{r\circ}\cdot v)$ which is possible only if $v=1$ and $G^{r\circ}=G^{s\circ}$. We necessarily have $J_G(rv)=J_G(s)$, whence the statement.

\noindent (5) Let $S$ be the sheet associated with $(M,\,Z^\circ t, M\cdot u)$. The elements in the Jordan class corresponding to this triple have stabilizer conjugate to $M^u$. Hence, if $S\subset G_{(n)}$ then $n=\dim G-\dim M^u$.

Moreover, $\overline{S}=\overline{G\cdot (Z(M)^{reg}_tu)}$ so 
$\dim S=\dim G+\dim Z(M)^{reg}_tu-\dim \sigma^{-1}(zu)$, where $\sigma$ is the natural map $\sigma\colon G\times Z(M)^{reg}_tu\to G\cdot (Z(M)^{reg}_tu)$ and $zu$ is a suitable element in $Z(M)^{reg}_tu$. We have
$$\begin{array}{rl}
\sigma^{-1}(zu)&=\{(g,ru)\in G\times Z(M)^{reg}_tu~|~ g\cdot ru=zu\}\\
&=\{(g,ru)\in G^u\times Z(M)^{reg}_tu~|~ g\cdot r=z\}\\
&=\bigcup_{i}\{(gg_i,\,z_iu)\in G^{zu}g_i\times(Z(M)^{reg}_t\cap G^u\cdot z)~|~g_i\cdot z_i=z\}\\
&=\bigcup_i (G^{su}g_i,\,z_iu)
\end{array}$$ where the union is over  finitely many elements because conjugate elements in a maximal torus lie in the same $W$-orbit.
Hence,  $\dim \sigma^{-1}(zu)=\dim G^{su}=\dim M^u$. 
 \hfill$\Box$

\bigskip

\begin{corollary}Let $H$ be a semisimple group. A conjugacy class $H\cdot su$ is rigid if and only if $s$ is exceptional and $G^{s\circ}\cdot u$ is rigid in $G^{s\circ}$. 
\end{corollary}
\pf
The Jordan class of $su$ is maximal if and only if $G^{s\circ}\cdot u$ is rigid in $G^{s\circ}$. Moreover, the corresponding sheet $S\subset G_{(n)}$ coincides with $G\cdot su$ if and only if $\dim S=n$. By Theorem~\ref{triplets} (e),  this happens if and only if $\dim Z(G^{s\circ})^\circ=0$. \hfill$\Box$

Equivalently we may state the classification theorem in the following form, which resembles the Lie algebra case. We restrict to the semisimple case, the general reductive case easily follows.

\begin{theorem}\label{pairs}Let $H$ be a semisimple group over an algebraically closed field of good characteristic. Then, 
\begin{enumerate}
\item[(1)] The sheets of $H$ are in one-to-one correspondence with $H$-orbits of pairs $(L, \OO)$
where $L$ is a Levi subgroup of $H$ and $\OO$ is a rigid conjugacy class in $[L,\,L]$.
\item[(2)] The sheets containing a unipotent conjugacy class up to a central factor are precisely those for which $\OO$ is unipotent up to a central factor. The sheets containing a genuine unipotent conjugacy class are those for which $\OO$ is unipotent.
\item[(3)] If a sheet contains a unipotent conjugacy class up to a central element, the involved unipotent class is unique.
\item[(4)]  The sheets containing a semisimple element are those for which $\OO$ is semisimple and every semisimple element lies in a unique sheet.
\item[(5)] The sheet associated with $(L,\,\OO)$ lies in $H_{(n)}$ with $n=\dim H-\dim L+\dim\OO$. 
Its dimension is equal to $n+\dim Z(L)^\circ$. 
\end{enumerate}
\end{theorem}
\pf It is straightforward to check that the correspondence between triples $(M, Z^\circ t, M\cdot u)$ as in Theorem~\ref{triplets} and pairs $(L,\,[L,\,L]\cdot su)$ as in Theorem~\ref{pairs}  where $L=C_G(Z(M)^\circ)$ and $s\in (Z^\circ t)^{reg}\cap [L,\,L]$ is a bijection on $G$-orbits.  \hfill$\Box$

 When $G=PGL_n({\mathbb C})$ every pseudo-Levi is a Levi subgroup and its center is always connected, so the parametrization in Theorem \ref{triplets} coincides with the parametrization of sheets of ${\mathfrak sl}_n({\mathbb C})$ as in \cite{bo}, and dimensions on corresponding sheets coincide.  For $G$ not containing only components of type $A_n$, there always exist sheets that do not have a counterpart in the corresponding Lie algebra because, among others, all  exceptional semisimple conjugacy classes are sheets.

\begin{remark}Theorem \ref{triplets} indicates how to classify sheets in a simple group $G$. Indeed, conjugacy classes of pseudo-Levi subgroups may be classified by using Remark \ref{deriziotis} and \cite[\S 2.2]{sommers}, the classes $Z^\circ t$ satisfying the condition on the centralizers  may be classified by using \cite[Proposition 15]{mcninch-sommers}, rigid unipotent elements can be classified by using the existing classification of rigid nilpotent orbits in \cite{el, ke, sp}. This, together with the computation of the dimension of each sheet and of each $G_{(n)}$, will be the content of a forthcoming paper. 
\end{remark}

\section{Acknowledgements}
The authors wish to express their gratitude to Mauro Costantini for  pointing out the content and the reference in Remark \ref{slodowy} and to Andrea Maffei for useful comments and suggestions.
The first named author was partially supported by Project CPDA071244 of the University of Padova.

After acceptance and typesetting of the paper, we discovered that the notions of Levi envelope, Jordan class and exceptional element in an algebraic group had already appeared in the paper \cite{lus} of G. Lusztig.

\end{document}